\documentclass[a4paper,twocolumn]{article} 

\usepackage[utf8]{inputenc}
\usepackage{amssymb,amsmath,amsthm} % math
\usepackage{caption, subcaption} % caption
\usepackage{epsfig} % for postscript graphics files
\usepackage[numbers]{natbib} % bib
\usepackage{url} % bib url
\newtheorem{prop}{Proposition} % prop
\newtheorem{lem}[prop]{Lemma} % lem
\newtheorem{exmp}[prop]{Example} % exmp
\usepackage[norelsize,algo2e,linesnumbered,ruled,vlined]{algorithm2e} % algorithms

\title{State space sets with common optimal feedback laws for nonlinear MPC} 
\author{Ruth Mitze, Raphael Dyrska, Kai K\"{o}nig and Martin M\"onnigmann\\
	Automatic Control and Systems Theory, Dep. of Mechanical Engineering,\\
	Ruhr-Universit\"at Bochum, 44801 Bochum, Germany.\\ E-mail: {\tt\small ruth.mitze@rub.de}, {\tt\small raphael.dyrska@rub.de}, {\tt\small kai.koenig-h4d@rub.de} and {\tt\small martin.moennigmann@rub.de}}
	
\begin{document}
\maketitle

\begin{abstract}                % Abstract of not more than 250 words.
In model predictive control (MPC), an optimal control problem (OCP) is solved for the current state and the first input of the solution, the optimal feedback law, is applied to the system. This procedure requires to solve the OCP in every time step.
Recently, a new approach was suggested for linear MPC. The parametric solution of a linear quadratic OCP is a piecewise-affine feedback law. The solution at a point in state space provides an optimal feedback law and a domain on which this law is the optimal solution. As long as the system remains in the domain, the law can be reused and the calculation of an OCP is avoided. In some domains the optimal feedback laws are identical. By uniting the corresponding domains, bigger domains are achieved and the optimal feedback law can be reused more often.
In the present paper, we investigate in how far this approach can be extended from linear to nonlinear MPC, we propose an algorithm and we illustrate the achieved savings with an example. \\

\end{abstract}

%===============================================================================

\section{Introduction}

% intro MPC
The ability of model predictive control (MPC) to consider constraints directly in the formulation makes MPC a favorable control scheme for many systems. However, MPC is computationally demanding because the optimal feedback law is calculated by periodically solving an optimal control problem (OCP) on a receding horizon. An entire field of research focusses on the reduction of the computational effort.

% Regional MPC
Regional MPC approaches intend to reduce the computational effort by reducing the number of OCPs that are solved.
Most publications focus on the linear case~\cite{Jost2015a,Konig2017,Konig2018,Berner2018} and exploit characteristics of the parametric solution of the OCP which is a piecewise-affine feedback law~\cite{Bemporad2002,Seron2003}. It is the central idea of regional MPC to reuse the optimal feedback law of the previous state if the current state is located in the same region as the previous state.  In contrast to explicit approaches~\cite{Bemporad2002,Tondel2003,Johansen2002}, regional MPC approaches do not require the parametric solution itself.
It has been shown in~\cite{Monnigmann2015} that regional MPC approaches can be extended to the nonlinear case. 

% Kais approach
Recently, a regional MPC approach was presented for linear MPC. 
The approach exploits that in some regions of the piecewise-affine feedback law, the same optimal feedback law applies.
It presents a simple criterion for detecting if the current state is part of a region with the same optimal feedback law as the previous state~\citep{Konig2020}. This enables reusing the optimal feedback law of the previous state more often and therefore to avoid solving an OCP more often. 
In contrast to other approaches \citep{Geyer2008,Kvasnica2012,Kvasnica2013}, it does not require an explicit solution.

% us
It is the purpose of the present paper to extend the approach from~\cite{Konig2020} from the linear to the nonlinear case.
We proceed analogously to the linear case and identify subsets of active sets that already define the optimal feedback law. Any optimal active set that contains such a subset will then identify the same optimal feedback law as the subset. 
Then, we present an approach that reduces the number of solved OCPs for nonlinear MPC. Due to the limitations that arise caused by the nonlinearities, the proposed approach is different from the approach from~\cite{Konig2020} and requires offline calculations.
Finally, the approach is illustrated with an example.

This paper is structured as follows. Section~\ref{sec:ProblemStatement} introduces the class of nonlinear optimal control problems treated here. Section~\ref{sec:commonFeedbackLaws} describes state space sets with common feedback laws. Implementational aspects and an example are discussed in Sects.~\ref{sec:implementation} and~\ref{sec:example}, respectively. Brief conclusions are given in Sect.~\ref{sec:conclusion}.

\subsection*{Notation}

% only particular lines of matrix
For any $M\in\mathbb{R}^{a\times b}$ and any ordered set $\mathcal{M}\subseteq\{1,...,a\}$ let $M_{\mathcal{M}}\in\mathbb{R}^{\vert\mathcal{M}\vert\times b}$ be the submatrix of $M$ containing all rows indicated by $\mathcal{M}$.

\section{Problem statement and preliminaries} \label{sec:ProblemStatement}

% system
Consider a discrete-time system
\begin{align} \label{eq:System}
x(k+1)=f(x(k),u(k)),\, k= 0, 1, \dots
\end{align}
that must respect constraints of the form
\begin{align*}
  u(k)\in\mathcal{U}\subseteq\mathbb{R}^{m}, \,
  x(k)\in\mathcal{X}\subseteq\mathbb{R}^{n}, \,
  k= 0, 1, \dots
\end{align*}
with input variables $u(k)\in\mathbb{R}^m$, state variables $x(k)\in\mathbb{R}^n$, and a nonlinear function $f:\mathbb{R}^n\times\mathbb{R}^m\rightarrow\mathbb{R}^n$. We assume $f$ is twice continuously differentiable, $f(0,0)=0$ holds, and $\mathcal{U}$ and $\mathcal{X}$ are compact full-dimensional polytopes that contain the origin in their interiors.

% OCP
The optimal control problem (OCP) treated in the present paper reads
\begin{align} \label{eq:OCP}
\begin{split}
\min_{U,X} \quad &x(N)^T Px(N)+\sum_{k=0}^{N-1}\left(x(k)^TQx(k)+u(k)^TRu(k)\right) \\
\textrm{s.t.} \quad 
&x(k+1)=f(x(k),u(k)), \: k=0,...,N-1\\
&u(k)\in\mathcal{U}, \: k=0,...,N-1\\
&x(k)\in\mathcal{X}, \: k=0,...,N-1\\
&x(N)\in\mathcal{T},
\end{split}
\end{align}
where $U=\left(u^T(0),...,u^T(N-1)\right)^T\in\mathbb{R}^{Nm}$ and $X=\left(x^T(1),...,x^T(N)\right)^T\in\mathbb{R}^{Nn}$ collect the inputs and states, respectively, the initial state $x(0)$ is given, $Q,P\in\mathbb{R}^{n\times n}$ and  $R\in\mathbb{R}^{m\times m}$, with $Q,P,R\succ 0$ are the usual weighting matrices, $N\in\mathbb{N}$ is the horizon, and $\mathcal{T}\subseteq\mathcal{X}$ is the terminal set. 

% constraint indexing
We assume the set $\mathcal{U}$ is defined by a finite number of halfspaces. Let $q_\mathcal{U}$ denote the number of halfspaces that bounds $\mathcal{U}$. Furthermore, let $q$ and $\mathcal{Q}=\{1,...,q\}$ refer to the total number of constraints in~\eqref{eq:OCP} and their index set, respectively.

% OCP after substitution
By substituting~\eqref{eq:System} into~\eqref{eq:OCP}, the OCP~\eqref{eq:OCP} can be stated in the form
\begin{align} \label{eq:newOCP}
\begin{split}
\min_{U} \quad &V(x(0),U) \\
\textrm{s.t.} \quad &G(x(0),U)\leq 0.
\end{split}
\end{align}

% KKT
The KKT-conditions that solve~\eqref{eq:newOCP} read
\begin{align}\label{eq:KKT}
\begin{split}
\nabla_U\left(V(x(0),U)+\sum_{i=1}^{q}\lambda_iG_i(x(0),U)\right)=0,\\
\lambda_iG_i(x(0),U)=0,\:i=1,...,q,\\
G_i(x(0),U)\leq 0,\:i=1,...,q,\\
\lambda_i\geq 0,\:i=1,...,q,
\end{split}
\end{align}
where $\lambda\in\mathbb{R}^q$ denotes the Lagrangian multipliers. Let $\mathcal{F}$ refer to the set of initial states $x(0)$ such that~\eqref{eq:KKT} has a solution. For any $x(0)\in\mathcal{F}$, let $U^\star(x(0))$, $\lambda^\star(x(0))$ denote the optimal solution to~\eqref{eq:KKT}, $U^\star:\mathcal{F}\rightarrow\mathbb{R}^{Nm}$, $\lambda^\star:\mathcal{F}\rightarrow\mathbb{R}^{q}$. We call $U^\star(x(0))$ the \textit{optimal control law}.

% optimal feedback law
In model predictive control, problem~\eqref{eq:newOCP} is solved in every time step for the current state $x(0)$ and the first input of the optimal control law $U_{\{1,...,m\}}^\star(x(0))$ is applied to system~\eqref{eq:System}. We call the first input of the optimal control law the \textit{optimal feedback law}
\begin{align*}
u^\star(x(0)):=U_{\{1,...,m\}}^\star(x(0)),
\end{align*}
where $u^\star: \mathcal{F}\rightarrow\mathbb{R}^{m}$.

We state an important fact about the solution structure in Proposition~\ref{prop:solutionStructure} below. As preparation, we introduce the following sets.

% active and inactive set
For any $x(0)\in\mathcal{F}$, let $\mathcal{A}(x(0))$ and $\mathcal{I}(x(0))$ refer to the optimal active set $\mathcal{A}(x(0))=\left\{i\in\mathcal{Q}\vert G(x(0),U^\star(x(0)))=0\right\}$ and the corresponding inactive set $\mathcal{I}(x(0))=\mathcal{Q}\backslash\mathcal{A}(x(0))$. 
% weakly and strongly active
The sets $\mathcal{W}(x(0))$ and $\mathcal{S}(x(0))$ refer to the sets of weakly active constraints $\mathcal{W}(x(0))=\left\{i\in\mathcal{A}(x(0))\,\vert\,\lambda_i^\star(x(0))=0\right\}$ and strongly active constraints $\mathcal{S}(x(0))=\mathcal{A}(x(0))\backslash\mathcal{W}(x(0))$, respectively.
We often drop the argument $x(0)$ for brevity. 

% active set exists
We say an active set $\mathcal{A}$ exists for the problem~\eqref{eq:newOCP} if $\mathcal{A}$ appears as the active set for an $x(0)\in\mathcal{F}$. 

% KKT with sets
With the sets just introduced,~\eqref{eq:KKT} can be expressed as
\begin{align}\label{eq:KKTwithSets}
\begin{split}
\nabla_U\left(V(x(0),U)+\sum_{i\in\mathcal{A}}\lambda_iG_i(x(0),U)\right)=0,\\
G_\mathcal{A}(x(0),U)=0,\\
\lambda_{\mathcal{I}\cup\mathcal{W}}=0,\\
G_\mathcal{I}(x(0),U)<0,\\
-\lambda_\mathcal{S}<0.
\end{split}
\end{align}
We denote the left-hand side of the equality and inequality constraints in~\eqref{eq:KKTwithSets} with $F_{\rm eq}^\mathcal{A}(U,\lambda,x(0))$ and $F_{\rm ineq}^{\mathcal{A},\mathcal{W}}(U,\lambda,x(0))$, respectively.

% active set --> U and region
\begin{prop}[{\cite[Prop.~1, Lem.~2]{Monnigmann2015}}]\label{prop:solutionStructure}
Consider an active set $\mathcal{A}$ such that a solution to $F_{\rm eq}^\mathcal{A}(U,\lambda,x(0))=0$ exists. Let the solution be denoted $(U_{\rm sol},\lambda_{\rm sol},x(0)_{\rm sol})$.

If $\frac{\partial F_{\rm eq}^\mathcal{A}(U,\lambda,x(0))}{\partial (U^T,\lambda^T)^T}\vert_{\rm sol}$ has full rank,
%\begin{itemize}
%\item $\mathcal{W}(x(0)_{\rm sol})=\emptyset$,
%\item the column vectors $\nabla_UG_\mathcal{A}(U_{\rm sol},x(0)_{\rm sol})$ are linear independent and
%\item ${\rm im}(\nabla_U^2\left(V(U_{\rm sol},x(0)_{\rm sol})-\sum_{i=1}^{q}(\lambda_{\rm sol})_iG_i(U_{\rm sol},x(0)_{\rm sol})\right)$ and ${\rm span}(\nabla_UG_\mathcal{A}(U_{\rm sol},x(0)_{\rm sol}))$ are orthogonal to one another,
%\end{itemize}
then $F_{\rm eq}^\mathcal{A}(U,\lambda,x(0))=0$ implicitly defines $U^\star(x(0))$ and $\lambda^\star(x(0))$ for those $x(0)$ on a region $\mathcal{R}(\mathcal{A})$,
\begin{align*}
\mathcal{R}(\mathcal{A})=\left\{x(0)\in\mathbb{R}^n\:\vert\:F_{\rm ineq}^{\mathcal{A},\emptyset}(U,\lambda,x(0))<0\right\}.
\end{align*}
\end{prop}

It follows from Prop.~\ref{prop:solutionStructure} that an active set that satisfies the conditions in Prop.~\ref{prop:solutionStructure} implicitly defines the optimal control law $U^\star$ for all states $x(0)\in\mathcal{R}(\mathcal{A})$.

\section{State space sets with common optimal feedback laws}\label{sec:commonFeedbackLaws}

% intro
In this section we derive a condition under that a subset of an active set defines the optimal feedback law $u^\star$. It follows that every active set with the same subset defines the same optimal feedback law.

% constraint order
We assume the constraints in~\eqref{eq:OCP} are ordered such that the input constraints on $u(0)$ appear first,
\begin{align} \label{eq:order}
\left.\begin{array}{c}
  u(0)\in\mathcal{U},\\ 
  \vdots 
\end{array}\right. .
\end{align}
Let the constraint order~\eqref{eq:order} be preserved in the constraints in~\eqref{eq:newOCP}. 

% constraints that depend on u(0)
\begin{lem}\label{lem:u0constraints}
Consider~\eqref{eq:newOCP} and assume the constraints are ordered as in~\eqref{eq:order}. Then, the constraints from~\eqref{eq:newOCP} can be stated in the form
\begin{align}\label{eq:splitG}
\left(\begin{array}{c}
  \tilde{G}\cdot u(0)-\tilde{w}\\
  G_{\{q_\mathcal{U}+1,...,q\}}(x(0),U)
\end{array}\right)
  \leq 0,
\end{align}
with $\tilde{G}\in\mathbb{R}^{q_\mathcal{U}\times m}$ and $\tilde{w}\in\mathbb{R}^{q_\mathcal{U}}$.

\begin{proof}
The first constraint in~\eqref{eq:order} is
\begin{align*}
  u(0)\in\mathcal{U}.
\end{align*}
By assumption, $\mathcal{U}$ is a polytope and bounded by $q_\mathcal{U}$ halfspaces. Therefore, the first $q_\mathcal{U}$ rows of $G(x(0),U)$ are linear and only depend on $u(0)$. They can be expressed in the form
\begin{align*}
\tilde{G}\cdot u(0)\leq\tilde{w}.
\end{align*}
\hfill$\square$
\end{proof}
\end{lem}

% tilde A --> u(0)
Consider an arbitrary active set $\mathcal{A}$ that exists for~\eqref{eq:newOCP}. Let $\tilde{\mathcal{A}}$ contain the indices that are active in the first rows of~\eqref{eq:splitG},
\begin{align}\label{eq:tildeA}
\tilde{\mathcal{A}}:=\mathcal{A}\cap\{1,...,q_\mathcal{U}\}.
\end{align}
With Lem.~\ref{lem:u0constraints}, the constraints that correspond to the active set $\tilde{\mathcal{A}}$ are
\begin{align}\label{eq:tildeConstraints}
\tilde{G}_{\tilde{\mathcal{A}}}\cdot u(0)-\tilde{w}_{\tilde{\mathcal{A}}}=0.
\end{align}
If~\eqref{eq:tildeConstraints} already determines the optimal feedback law $u(0)$, then this is the optimal feedback law for all active sets $\mathcal{A}^\prime$ that exist for~\eqref{eq:newOCP} and satisfy
\begin{align*}
\mathcal{A}^\prime\cap\{1,...,q_\mathcal{U}\}=\tilde{\mathcal{A}}.
\end{align*}
We summarize these active sets in the set
\begin{align*}
\mathcal{M}(\mathcal{A})=\left\{\mathcal{A}^\prime\subseteq\mathcal{Q}\,\vert\,\mathcal{A}^\prime\cap\{1,...,q_\mathcal{U}\}=\mathcal{A}\cap\{1,...,q_\mathcal{U}\}\right\}. 
\end{align*}
It follows that the same optimal feedback law applies to the union of regions
\begin{align}\label{eq:regionUnion}
\Gamma(\mathcal{A})=\bigcup_{\mathcal{A}^\prime\in\mathcal{M}(\mathcal{A})}\mathcal{R}(\mathcal{A}^\prime).
\end{align}
We point out that this equality only holds for the first entry of the optimal control law, the optimal feedback law $u(0)$. Of course, the control law $u(0),...,u(N-1)$ might differ for each region $\mathcal{R}(\mathcal{A})$ in~\eqref{eq:regionUnion}. We summarize the explanations given so far in the following proposition.

% same u(0), linear part
\begin{prop}\label{prop:tildeIdenticalFeedbackLaw}
Consider an arbitrary active set $\mathcal{A}$ that exists for~\eqref{eq:newOCP} and assume the constraints are ordered as in~\eqref{eq:order}. Let $\tilde{\mathcal{A}}$, $\tilde{G}$ and $\tilde{w}$ be defined as in~\eqref{eq:tildeA} and~\eqref{lem:u0constraints}.

If $\vert\tilde{\mathcal{A}}\vert=m$ and if $\tilde{G}_{\tilde{\mathcal{A}}}$ has full rank, then 
\begin{align}\label{eq:tildeControlLaw}
u^\star=\tilde{G}_{\tilde{\mathcal{A}}}^{-1}\cdot\tilde{w}_{\tilde{\mathcal{A}}}
\end{align}
is the optimal feedback law for all $x\in\Gamma(\mathcal{A})$, where $\Gamma(\mathcal{A})$ is as in~\eqref{eq:regionUnion}.

\begin{proof}
The constraints can be formulated as in~\eqref{eq:splitG}, since the assumptions for Lem.~\ref{lem:u0constraints} hold. All constraints $i\in\mathcal{A}$ hold with equality and $\tilde{\mathcal{A}}\subseteq\mathcal{A}$, therefore~\eqref{eq:tildeConstraints} holds, with $\tilde{G}_{\tilde{\mathcal{A}}}\in\mathbb{R}^{\vert\tilde{\mathcal{A}}\vert\times m}$ and $\tilde{w}_{\tilde{\mathcal{A}}}\in\mathbb{R}^{\vert\tilde{\mathcal{A}}\vert}$. 
By assumption, $\tilde{G}_{\tilde{\mathcal{A}}}$ has full rank and is a square matrix, because $\vert\tilde{\mathcal{A}}\vert=m$. Hence, the inverse $\tilde{G}_{\tilde{\mathcal{A}}}^{-1}$ exists.
Reformulating~\eqref{eq:tildeConstraints} results in
\begin{align}\label{eq:helperProof}
u(0)=\tilde{G}_{\tilde{\mathcal{A}}}^{-1}\cdot\tilde{w}_{\tilde{\mathcal{A}}}
\end{align}
which provides~\eqref{eq:tildeControlLaw}. It follows that under the stated assumptions the active constraints $\tilde{\mathcal{A}}$ define the optimal feedback law $u(0)$.

For all $x\in\Gamma(\mathcal{A})$ holds $\tilde{\mathcal{A}}\subseteq\mathcal{A}$~\eqref{eq:regionUnion} which implies~\eqref{eq:tildeConstraints}. Therefore, the optimal feedback law~\eqref{eq:helperProof} is the optimal feedback law for all $x\in\Gamma(\mathcal{A})$.\hfill$\square$
\end{proof}
\end{prop}

\section{Implementational aspects}\label{sec:implementation}

% intro
Proposition~\ref{prop:tildeIdenticalFeedbackLaw} may show that regions with common optimal feedback laws also exist for the nonlinear case. The criterion for that is similar to the criterion for the linear case~\cite[Prop.~1, Lem.~3]{Konig2020}. It is an obvious question whether the regional MPC approach for the linear case from~\cite[Sect.~4]{Konig2020} can be transferred to the nonlinear case similarly.

% calculating the regions
The approach presented in~\cite[Sect.~4]{Konig2020} determines regions in state space with common optimal feedback laws online. In the linear case, the regions are polytopes~\cite[Sect. 4.1]{Bemporad2002}. The computational effort to calculate the polytope defined by an active set is relatively small, it only comprises matrix operations (see e.g.~\cite[Lem.~2]{Jost2015a}).
In the nonlinear case, the regions are bounded nonlinearly. Without giving details we claim the computational effort to calculate the region defined by an active set is relatively high. Therefore, calculating those regions online is not an option. Instead, we suggest determining the regions with common optimal feedback laws offline. 
%Note that this does not require to calculate the entire explicit solution. 
This requires to calculate those regions that are defined by active sets $\mathcal{A}$ such that for their subset $\tilde{\mathcal{A}}$~\eqref{eq:tildeA} holds $\vert\tilde{\mathcal{A}}\vert=m$, and $\tilde{G}_{\tilde{\mathcal{A}}}$ has full row rank (see conditions in Prop.~\ref{prop:tildeIdenticalFeedbackLaw}). 

% underestimate with ellipsoid, save feedback law
Regions defined by active sets with the same subset $\tilde{\mathcal{A}}$ are then united and for each of those unions, the optimal feedback law $u^\star$ is determined with~\eqref{eq:tildeControlLaw}.
Further, each union is underestimated by one or more ellipsoids $\mathcal{E}$ of the form
\begin{align}\label{eq:ellipsoid}
\mathcal{E}=\left\{x\in\mathbb{R}^n\,\vert\, (x-x_C)^TE(x-x_C)\leq 1\right\},
\end{align} 
with $E\in\mathbb{R}^{n\times n}$ where $x_C\in\mathbb{R}^n$ denotes the center of the ellipsoid.
Finally, all pairs consisting of the ellipsoid and the corresponding optimal feedback law $(\mathcal{E},u^\star)$ are collected in the set $\mathcal{S}$.

% online algorithm
The online part of the approach is shown in Alg.~\ref{algorithm}. It exploits that the optimal feedback law $u^\star$ is known for all states that are an element of an ellipsoid $\mathcal{E}$. 
The algorithm tests if the current state is part of one of the ellipsoids (lines 2,3) and, if appropriate, the optimal feedback law is set to the optimal feedback law that corresponds to that specific ellipsoid (line 4). The OCP is only solved otherwise (line 6).
\begin{algorithm2e}[h]
\textbf{Input:} current state $x(0)$, $\mathcal{S}$ (determined offline) \\
\For{every $(\mathcal{E}_i,u^\star_i)\in\mathcal{S}$}{
	\If{$x(0)\in\mathcal{E}_i$}{
		set $u^\star(x(0))=u^\star_i$
	}
}
\If{$u^\star(x(0))$ is undefined}{
	solve~\eqref{eq:newOCP} for $u^\star(x(0))$
}
\textbf{Output:} $u^\star(x(0))$
\caption{NMPC using predefined regions with corresponding optimal feedback laws\label{algorithm}}
\end{algorithm2e}

% why underestimation
The fact that the united regions were underestimated by ellipsoids lowers the performance of the proposed algorithm because fewer states are identified to be part of a region with known optimal feedback law $u^\star$. The underestimation is necessary to keep the computational effort for the membership test in line 2 low, solving simple inequalities is sufficient to test whether a state is part of an ellipsoid~\eqref{eq:ellipsoid}. This is essential because in case none of the ellipsoids applies we still solve the OCP (line 6). Note that the underestimation does not affect the optimal feedback law which is the output of Alg.~\ref{algorithm}.

\section{Example}\label{sec:example}

% intro
In this section, we illustrate Prop.~\ref{prop:tildeIdenticalFeedbackLaw} and the approach presented in Sect.~\ref{sec:implementation} with the following example.

% example
\begin{exmp}\label{exmp:Pannocchia2011}
Consider the system~\cite[Sect. $\mathord{\mathrm{VI}}$A, C2]{Pannocchia2011}
\begin{align*}
x(k+1)=\left(\begin{matrix}x_1(k)+u(k)\\bx_2(k)+u(k)^3\end{matrix}\right),
\end{align*}
with $\vert u(k)\vert\leq 1$, $\mathcal{T}=\left\{ x(k)\in\mathbb{R}^2\vert x(k)^TPx(k)\leq\alpha\right\}$, $Q=I^2$, $R=1$ and $N=3$, where $I$ denotes the identity matrix, $b=0.9$, $P=\left(\begin{matrix}4&0\\0&10.53\end{matrix}\right)$ and $\alpha=1.1$.
The feasible set for this example is shown in Fig.~\ref{fig:activeSet}. It was generated by determining the solution for a finite number of states $x(0)$ and consists of $26$ different optimal active sets.
\begin{figure}[tbh]
    \begin{center}
	\includegraphics[trim=0px 0 0px 0, clip, width=.45\textwidth]{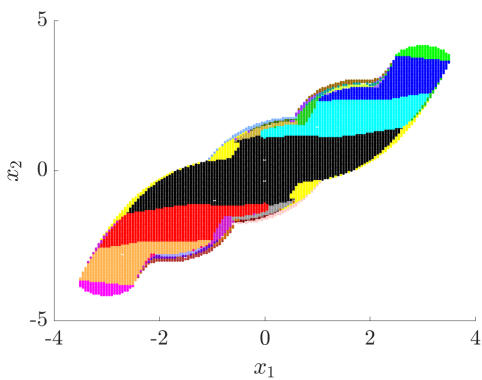}
    \caption{Feasible set for Example~\ref{exmp:Pannocchia2011}. Different colors indicate states with different optimal active sets.}
    \label{fig:activeSet}
    \end{center}
\end{figure}
\end{exmp}

\begin{figure}[tbh]
    \begin{center}
	\includegraphics[trim=0px 0 0px 0, clip, width=.45\textwidth]{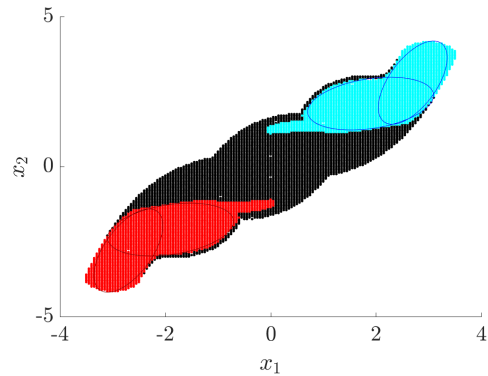}
    \caption{Feasible set for Example~\ref{exmp:Pannocchia2011}. For states that are illustrated by blue and red color, the optimal feedback law is $u^\star=-1$ and $u^\star=1$, respectively. Ellipsoids underestimate blue and red regions.}
    \label{fig:saturated}
    \end{center}
\end{figure}

% illustration of Proposition
For Example~\ref{exmp:Pannocchia2011} the sets $\tilde{\mathcal{A}}$ such that the conditions in Prop.~\ref{prop:tildeIdenticalFeedbackLaw} hold result $\tilde{\mathcal{A}}=\{1\}$ and $\tilde{\mathcal{A}}=\{2\}$. The corresponding optimal feedback laws with~\eqref{eq:tildeControlLaw} are $u^\star=-1$ and $u^\star=1$, respectively. 
In this example, seven different optimal active sets in the solution contain the subset $\{1\}$ and seven contain the subset $\{2\}$. Fig.~\ref{fig:saturated} shows all states such that the optimal active set contains the subsets $\{1\}$ and $\{2\}$ by blue and red color, respectively, states such that the optimal active set does not contain any of these subsets are shown in black. With Prop.~\ref{prop:tildeIdenticalFeedbackLaw}, all optimal active sets that contain the subsets $\{1\}$ and $\{2\}$ define the optimal feedback laws $u^\star=-1$ and $u^\star=1$, respectively.

% reduction of computational effort
The red and blue regions in Fig.~\ref{fig:saturated} are underestimated with two ellipsoids each. The proportion of the ellipsoids on the feasible set is $43.7$\%. 
Thus, for $43.7$\% of the states $x(0)$, Alg.~\ref{algorithm} sets the optimal feedback law to a predefined value and no OCP is solved.

\begin{figure}[tbh]
    \begin{center}
	\includegraphics[trim=0px 0 0px 0, clip, width=.45\textwidth]{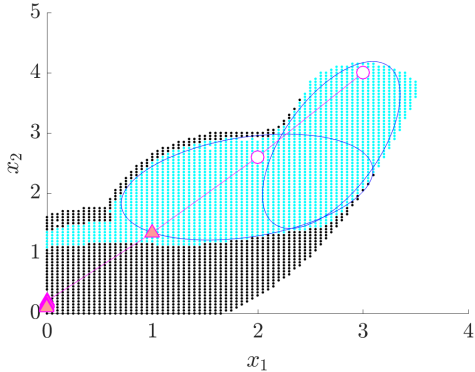}
    \caption{Detail of Fig.~\ref{fig:saturated} with closed-loop trajectory that results for $x(0)=(3,4)^T$ for Example~\ref{exmp:Pannocchia2011}. Trajectory states that are part of an ellipsoid are marked with white circles, all other states are marked with light red triangles.}
    \label{fig:trajectory}
    \end{center}
\end{figure}

% illustration of algorithm with sample trajectory
Fig.~\ref{fig:trajectory} illustrates a sample closed-loop trajectory. 
For all states that are included in an ellipsoid (white circles), the optimal feedback law is known and no OCP is solved (lines 3, 4 in Alg.~\ref{algorithm}). For all states that are not part of an ellipsoid (light red triangles), an OCP is solved to determine the optimal feedback law (lines 5, 6 in Alg.~\ref{algorithm}). 
Note that the third state of the trajectory is part of the blue region where the optimal feedback law is known, but it is not part of an ellipsoid. Therefore, an OCP is solved for this state.

\section{Conclusion}\label{sec:conclusion}

We presented a criterion for detecting when a subset of the optimal active set defines the optimal feedback law. It follows that active sets that contain the same subset define the same optimal feedback law. 
We used the criterion and proposed an algorithm to reduce the number of OCPs that are solved. The online reduction was illustrated with an example, the calculation of an OCP was avoided for more than $40$\% of the states.
A drawback of the approach presented in this paper is its required offline effort. The approach requires to calculate, unite, and underestimate regions with common optimal feedback laws offline. This limits the approach to low complexity problems.

\section{Acknowledgements}                           % Place acknowledgements
Support by the Deutsche Forschungsgemeinschaft (DFG) under grant MO 1086/15-1 is gratefully acknowledged.

\bibliographystyle{plain} 
\bibliography{root}             % bib file to produce the bibliography
                                                     % with bibtex (preferred)
                                                   
%\appendix
%\section{A summary of Latin grammar}    % Each appendix must have a short title.

\end{document}